\documentclass[11pt]{amsart}

\newlength{\basicwidth}
\newlength{\shortbasicwidth}
\newlength{\basicheight}
\numberwithin{equation}{section}

\begin{document}

\title{
{\bf{On a sine polynomial of Tur\'an}}}
\maketitle

\begin{center}
HORST ALZER$^a$  \, and  \, MAN KAM KWONG$^b$
\footnote{The research of this author is supported by the Hong Kong Government GRF Grant PolyU 5003/12P and the Hong Kong Polytechnic University Grants G-UC22 and G-UA10}
\end{center}

\vspace{1.5cm}
\begin{center}
$^a$
Morsbacher Str. 10, 51545 Waldbr\"ol, Germany\\
\emph{email:} \tt{h.alzer@gmx.de}
\end{center}

\vspace{0.4cm}
\begin{center}
$^b$
Department of Applied Mathematics, The Hong Kong Polytechnic University,\\
Hunghom, Hong Kong\\
\emph{email:} \tt{mankwong@polyu.edu.hk}
\end{center}

\vspace{2cm}
\noindent
{\bf{Abstract.}}
In 1935, P. Tur\'an proved that
$$
S_{n,a}(x)=
\sum_{j=1}^n{n+a-j\choose n-j} \sin(jx)>0 \quad{(n,a\in\mathbf{N}; 0<x<\pi).}
$$
We present various related inequalities. 
Among others, we show that the refinements
$$
S_{2n-1,a}(x)\geq \sin(x)
\quad\mbox{and}
\quad{S_{2n,a}(x)\geq 2\sin(x)(1+\cos(x))}
$$
are valid for all integers $n\geq 1$ and real numbers $a\geq 1$
 and $x\in(0,\pi)$. Moreover, we apply our theorems on sine sums to obtain inequalities for the Chebyshev polynomials of the second kind.

\vspace{1cm}

\vspace{0.3cm}
\newpage

\section{Introduction}

We define the sequence $\sigma_{n,k}(z)$ recursively by
$$
\sigma_{n,0}(z)=\sum_{j=0}^n z^j,
\quad{\sigma_{n,k}(z)=\sum_{j=0}^n \sigma_{j,k-1}(z)}
\quad{(k\in\mathbf{N})}.
$$
Then we have
$$
\sigma_{n,k}(e^{ix})=
\sum_{j=0}^n
{n+k-j\choose k}\cos(jx)+ i
\sum_{j=1}^n
{n+k-j\choose k}\sin(jx).
$$
In 1935, Tur\'an \cite{T1} studied the imaginary part and proved by induction on $n$ and $k$
 the remarkable inequality
\begin{equation}
\sum_{j=1}^n{n+k-j\choose k} \sin(jx)>0
\quad{(n,k\in\mathbf{N}; \, 0<x<\pi)}.
\end{equation}
In the same paper, Tur\'an presented
an elegant inequality for a sine sum in two variables. He demonstrated that the following companion of (1.1) is valid:
\begin{equation}
\sum_{j=1}^n{n+k-j\choose k}\frac{ \sin(jx)\sin(jy)}{j}>0
\quad{(n,k\in\mathbf{N}; \, 0<x,y<\pi)}.
\end{equation}

 Szeg\"o \cite{S} applied (1.1) with $k=2$ to establish a theorem on univalent functions.
An extension of (1.1) with $k=2$
was given by Alzer and Kwong \cite{AK}, whereas  Alzer and Fuglede \cite{AF}
  offered a positive lower bound for the sine polynomial
  in (1.1) under the assumption that $n,k\geq 2$. In fact, they
   proved that
\begin{equation}
\sum_{j=1}^n{n+k-j\choose k} \sin(jx)>\frac{x(\pi-x)}{\pi}
\quad{(2\leq n,k\in\mathbf{N}; \, 0<x<\pi)}.
\end{equation}

The definition of the sums given in (1.1) and (1.2) requires that $k$ is a nonnegative integer. However, the identity
$$
{n+k-j\choose k}
={n+k-j\choose n-j}
$$
reveals that if we use the second binomial coefficient then $k$ can be any real number. Therefore it is natural  to ask for all real parameters $a$ and $b$ such that we have for all integers $n\geq 1$ and real numbers $x\in (0,\pi)$,  $y\in (0,\pi)$:
\begin{equation}
S_{n,a}(x)=\sum_{j=1}^n{n+a-j\choose n-j} \sin(jx)>0
\quad{\mbox{and}}
\quad{\Theta_{n,b}(x,y)=\sum_{j=1}^n{n+b-j\choose n-j} \frac{\sin(jx)\sin(jy)}{j}>0.
}
\end{equation}
In this paper we solve both problems.  Moreover, we provide several  closely related inequalities. Among others, we
show that there exist functions $\lambda(x)$, $\mu(x)$ and $\lambda^*(x,y)$, $\mu^*(x,y)$ such that the estimates
$$
S_{2n-1,a}(x)\geq \lambda(x)>0,
\quad{S_{2n,a}(x)\geq \mu(x)>0}
$$
and
$$
\Theta_{2n-1,a}(x,y)\geq \lambda^*(x,y)>0,
\quad{\Theta_{2n,a}(x,y)\geq \mu^*(x,y)>0}
$$
are valid for all integers $n\geq 1$ and real numbers $a\geq 1$, $x\in (0,\pi)$, $y\in (0,\pi)$.

In the next section we collect some lemmas
which we need to prove  our main results given in Section 3. Finally, in Section 4 we apply our theorems  to obtain inequalities for sums involving Chebyshev polynomials of the second kind and we also offer  new integral inequalities for
these polynomials.
Throughout, we maintain the notations introduced in this section.

For more information on inequalities for trigonometric sums and polynomials we refer to
the monograph Milovanovi\'c, Mitrinovi\'c, Rassias \cite[chapter 6]{MMR}.

\vspace{0.5cm}
\section{Lemmas}

The first lemma is due to Fej\'er 
\cite{F}.

\vspace{0.5cm}
{\bf{Lemma 2.1.}} \emph{Let $x\in (0,\pi)$ and}
\begin{equation}
\phi_n(x)=2\sum_{j=1}^{n-1}  \sin(jx)+\sin(nx).
\end{equation}
\emph{Then, $\phi_n(x)>0$ for $n=1,2$ and $\phi_n(x)\geq 0$ for $n\geq 3$.}

\vspace{0.5cm}
To prove (2.1) Fej\'er made use of the identity
$$
\phi_n(x)=\sin(x)\Bigl(n+2\sum_{j=1}^{n-1}\sum_{k=1}^j\cos(kx)\Bigr)
$$
which he obtained by comparing the coefficients of certain power series. Here, we offer a different proof which is  more elementary than Fej\'er's approach.

\vspace{0.3cm}
\begin{proof}
Multiplying both sides of (2.1) by $\sin(x/2)$ gives
$$
\sin(x/2)\phi_n(x)=\sum_{j=1}^{n-1}2\sin(x/2)\sin(jx)+\sin(x/2)\sin(nx)
$$
$$
=\sum_{j=1}^{n-1}\bigl(\cos((j-1/2)x)-\cos((j+1/2)x)\bigr)+\frac{1}{2}\bigl(\cos((n-1/2)x)-\cos((n+1/2)x)\bigr)
$$
$$
=\cos(x/2)-\frac{1}{2}\cos((n-1/2)x)
-\frac{1}{2}\cos((n+1/2)x)
$$
$$
=\cos(x/2)-\cos(x/2)\cos(nx)=\cos(x/2)\bigl(1-\cos(nx)\bigr).
$$
Since $x\in (0,\pi)$, we conclude that $\phi_n(x)>0$ for $n=1,2$ and $\phi_n(x)\geq 0$ for $n\geq 3$.
\end{proof}

\vspace{0.5cm}
Next, we present an identity which will be a helpful tool not only to establish our inequalities for trigonometric sums but also to cover all cases of equality.

\vspace{0.5cm}
{\bf{Lemma 2.2.}} \emph{Let $c_k$ $(k=1,...,n)$ be real numbers and}
$$
\gamma_{k,n}=c_k+2\sum_{j=1}^{n-k} (-1)^j c_{j+k}
\quad{(k=1,...,n)}.
$$
\emph{Then,}
\begin{equation}
\sum_{j=1}^n c_j \sin(jx)=\sum_{j=1}^n
\gamma_{j,n}\phi_j(x),
\end{equation}
\emph{where $\phi_j(x)$ is defined in} (2.1).

\vspace{0.3cm}
\begin{proof}
We have
$$
\sum_{j=1}^n
\gamma_{j,n}\phi_j(x)=
\sum_{j=1}^n
c_j\phi_j(x)
+2 \sum_{j=1}^n \Bigl(\phi_j(x)
\sum_{k=1}^{n-j} (-1)^k c_{k+j}\Bigr)
$$
$$
=
\sum_{j=1}^n
c_j\phi_j(x)
+2 \sum_{j=2}^n \Bigl(c_j
\sum_{k=1}^{j-1} (-1)^{j-k} \phi_{k}(x)\Bigr)
=c_1\phi_1(x)
+\sum_{j=2}^n c_j \Bigl(\phi_j(x)+
2 \sum_{k=1}^{j-1} (-1)^{j-k} \phi_{k}(x)\Bigr)
$$
$$
=
c_1\sin(x)+\sum_{j=2}^n c_j \sin(jx)=\sum_{j=1}^n c_j\sin(jx).
$$
\end{proof}

\vspace{0.5cm}
{\bf{Remark 2.3.}} An application of Lemmas 2.1 and 2.2 leads to a result of Steinig \cite{ST}, who proved that the sine polynomial in (2.2) is nonnegative on $(0,\pi)$, if
$\gamma_{k,n}\geq 0$ $(k=1,...,n)$.

\vspace{0.5cm}
{\bf{Lemma 2.4.}} \emph{Let $a_k$ $(k=1,...,n; n\geq 3)$ be real numbers such that}
$$
2a_k\leq a_{k-1}+a_{k+1}
\quad{(k=2,...,n-1)}
\quad{and}
\quad{0\leq 2a_n\leq a_{n-1}.}
$$
\emph{Then, for} $x\in (0,\pi)$,
\begin{equation}
L_n(x)=
\sum_{j=1}^n a_j \sin(jx)-\sum_{j=1}^{n-2} a_{j+2} \sin(jx)\geq 0.
\end{equation}

\vspace{0.3cm}
\begin{proof}
Let $x\in (0,\pi)$. We define
$$
\tilde{c}_j=a_j -a_{j+2}
\quad{(j=1,...,n-2)},
\quad{\tilde{c}_{n-1} =a_{n-1}}, 
\quad{\tilde{c}_n =a_n.}
$$
Then we have
\begin{equation}
L_n(x)=\sum_{j=1}^n \tilde{c}_j \sin(jx).
\end{equation}
Let
$$
\tilde{\gamma}_k =a_k-2a_{k+1}+a_{k+2}
\quad{(k=1,...,n-2)},
\quad{\tilde{\gamma}_{n-1}=a_{n-1}-2a_n},
\quad{\tilde{\gamma}_n =a_n}.
$$
By assumption,
\begin{equation}
\tilde{\gamma}_k\geq 0
\quad{(k=1,...,n)}.
\end{equation}
We have
$$
\tilde{\gamma}_k=
\tilde{c}_k +
2\sum_{j=1}^{n-k} 
(-1)^j \tilde{c}_{j+k}
\quad{(k=1,...,n)},
$$
so that Lemma 2.2 implies
\begin{equation}
\sum_{j=1}^n \tilde{c}_j \sin(jx)=\sum_{j=1}^n\tilde{\gamma}_j \phi_j(x).
\end{equation}
Applying Lemma 2.1 and (2.5) reveals that the sum on the right-hand side of (2.6) is nonnegative. From (2.4) and (2.6) we obtain (2.3).
\end{proof}

\vspace{0.5cm}
{\bf{Remark 2.5.}} We assume that  $L_n(x_0)=0$ with $x_0\in (0,\pi)$. The   proof of Lemma 2.4 shows that

(i) if $n\geq 3$ and  $a_n>0$, then $\phi_n(x_0)=0$.

Moreover, since $\phi_2(x_0)>0$, we obtain that

(ii)  if $n=3$, then $a_2 -2 a_3=0$;

(iii)  if $n\geq 4$, then
{$a_2 -2 a_3 +a_4=0$}.

\vspace{0.5cm}
The following lemma plays an important role in the proofs of Theorems 3.1 and 3.2 given in the next section.

\vspace{0.5cm}
{\bf{Lemma 2.6.}} \emph{For all integers $n\geq 3$ and real numbers $a\geq 1$, $x\in(0,\pi)$ we have}
\begin{equation}
S_{n,a}(x)\geq S_{n-2,a}(x).
\end{equation}

\vspace{0.3cm}
\begin{proof}
Let $a\geq 1$. We define
$$
\tilde{a}_j={n+a-j\choose  n-j}
\quad{(j=1,...,n)}.
$$
Then,
$$
\tilde{a}_{j-1}-2 \tilde{a}_j+\tilde{a}_{j+1}
=\frac{a(a-1)}{(n-j+1)! }
\prod_{\nu=1}^{n-j-1} (n+a-j-\nu)\geq 0
\quad{(j=2,...,n-1),}
$$
$$
\tilde{a}_{n-1}-2\tilde{a}_n=a-1,
\quad{\tilde{a}_n=1}.
$$
Since
$$
S_{n,a}(x)-S_{n-2,a}(x)=\sum_{j=1}^n \tilde{a}_j \sin(jx)-\sum_{j=1}^{n-2} \tilde{a}_{j+2}\sin(jx),
$$
we conclude from Lemma 2.4 that (2.7) holds for $x\in (0,\pi)$.
\end{proof}

\vspace{0.5cm}
{\bf{Remark 2.7.}} Let $n\geq 3$, $a\geq 1$, $\tilde{x}_0\in (0,\pi)$ and $S_{n,a}(\tilde{x}_0)=S_{n-2,a}(\tilde{x}_0)$.  Remark 2.5 implies that $a=1$ and $\phi_n(\tilde{x}_0)=0$.

\vspace{0.5cm}
\section{Trigonometric sums}

Our first two theorems show that Tur\'an's inequality (1.1) can be refined if we assume that either $n$ is odd or $n$ is even.

\vspace{0.5cm}
{\bf{Theorem 3.1.}}
\emph{For all odd integers $n\geq 1$ and real numbers $a\geq 1$, $x\in (0,\pi)$, we have
\begin{equation}
S_{n,a}(x)\geq \sin(x).
\end{equation}
Equality holds if and only if $n=1$ or $n=3$, $a=1$, $x=2\pi/3$.}

\vspace{0.5cm}
\begin{proof}
 Let $n\geq 3$ be odd, $a\geq 1$ and $x\in (0,\pi)$. Applying Lemma 2.6 gives
$$
S_{n,a}(x)\geq S_{1,a}(x)=\sin(x).
$$
We discuss the cases of equality. A short calculation yields that
$$
S_{3,1}(2\pi/3)=\sin(2\pi/3)=\frac{1}{2}\sqrt{3}.
$$
We assume that
$$
S_{n,a}(x)=\sin(x)=S_{1,a}(x).
$$
Case 1. $n=3$.\\
Applying Remark 2.7 leads to $a=1$ and
$$
\phi_3(x)=\sin(x)(1+2\cos(x))^2=0.
$$
This gives $x=2\pi/3$.

Case 2. $n\geq 5$.\\
From Lemma 2.6 we conclude that
$$
S_{1,a}(x)=S_{n,a}(x)\geq S_{5,a}(x)\geq S_{3,a}(x)\geq S_{1,a}(x).
$$
Thus,
$$
S_{3,a}(x)=S_{1,a}(x)
\quad\mbox{and}
\quad{
S_{5,a}(x)=S_{3,a}(x).}
$$
Applying Remark 2.7 yields
$$
\phi_3(x)=0
\quad\mbox{and}
\quad{\phi_5(x)=0}.
$$
The first equation yields $x=2\pi/3$. But,
$
\phi_5(2\pi/3)=\sqrt{3}/2$.
A contradiction. It follows that equality holds in (3.1) if and only if $n=1$ or $n=3$, $a=1$, $x=2\pi/3$.
\end{proof}

\vspace{0.5cm}
{\bf{Theorem 3.2.}}
\emph{For all even integers $n\geq 1$ and real numbers $a\geq 1$, $x\in (0,\pi)$, we have
\begin{equation}
S_{n,a}(x)\geq 2 \sin(x)(1+\cos(x)).
\end{equation}
Equality holds if and only if $n=2$, $a=1$ or $n=4$, $a=1$, $x=\pi/2$.}

\vspace{0.5cm}
\begin{proof}
Let $a\geq 1$ and $x\in (0,\pi)$. Using (2.7) gives for even $n$,
\begin{equation}
S_{n,a}(x)\geq S_{2,a}(x)=(1+a)\sin(x)+\sin(2x)\geq 2\sin(x)(1+\cos(x)).
\end{equation}
We have
$$
S_{2,1}(x)=2\sin(x)(1+\cos(x))
\quad\mbox{and}
\quad{S_{4,1}(\pi/2)=2\sin(\pi/2)(1+\cos(\pi/2))=2.}
$$
Next, we assume that
\begin{equation}
S_{n,a}(x)=2\sin(x)(1+\cos(x)).
\end{equation}
Case 1. $n=2$.\\
We obtain
$$
0=S_{2,a}(x)-2\sin(x)(1+\cos(x))=(a-1)\sin (x).
$$
Thus, $a=1$.

Case 2. $n=4$.\\
Applying (3.3) and (3.4) yields 
$$
 S_{4,a}(x)=S_{2,1}(x).
$$
It follows from Remark 2.7 that $a=1$ and
$$
\phi_4(x)=8\sin(x)(1+\cos(x))\cos^2(x)=0.
$$
Thus, $x=\pi/2$.

Case 3. $n\geq 6$.\\
From (2.7), (3.3) and (3.4) we conclude that
$$
S_{4,a}(x)=S_{2,a}(x)
\quad\mbox{and}
\quad{S_{6,a}(x)=S_{4,a}(x),}
$$
so that Remark 2.7 gives
$$
\phi_4(x)=0
\quad\mbox{and}
\quad{\phi_6(x)=0.}
$$
The 
 first equation yields $x=\pi/2$, but $\phi_6(\pi/2)=2$. Hence, equality holds in (3.2) if and only if $n=2$, $a=1$ or $n=4$, $a=1$, $x=\pi/2$.
\end{proof}

\vspace{0.5cm}
{\bf{Remark 3.3.}} From Theorems 3.1 and 3.2 we conclude that the estimate
\begin{equation}
S_{n,a}(x)>\sin(x) \min\{1,2\,(1+\cos(x))\}
\end{equation}
is valid for $n\geq 1$, $a\geq 1$ and $x\in (0,\pi)$. The lower bounds given in (1.3) and (3.5) cannot be compared. Indeed, the function
$$
x\mapsto \frac{x(\pi-x)}{\pi}-
\sin(x)\min\{1,2\,(1+\cos(x))\}
$$
is negative on $(0,x_1)$ and positive on $(x_1,\pi)$, where $x_1=2.204...$.

\vspace{0.5cm}
The following extension of inequality (1.1) is valid.

\vspace{0.5cm}
{\bf{Theorem 3.4.}}
\emph{Let $a$ be a real number. The inequality 
\begin{equation}
S_{n,a}(x)>0
\end{equation}
holds for all integers $n\geq 1$ and real numbers $x\in (0,\pi)$
if and only if $a\geq 1$.}

\vspace{0.5cm}
\begin{proof}
From (3.5)  we conclude that if $a\geq 1$, then (3.6) holds for all $n\geq 1$ and $x\in (0,\pi)$. Conversely, let (3.6) be valid 
for all $n\geq 1$ and $x\in (0,\pi)$. Since $S_{2,a}(\pi)=0$, we obtain
$$
\frac{d}{dx} S_{2,a}(x)\Big{|_{x=\pi}}=1-a\leq 0.
$$
Thus, $a\geq 1$.
\end{proof}

\vspace{0.5cm}
Next, we present inequalities for the sine polynomial
$$
S^*_{n,a}(x)=\sum_{j=1 \atop j \,  odd}^n
{n+a-j\choose  n-j} \sin(jx).
$$
Applications   of Theorems 3.1 and 3.2 lead to counterparts of (3.1) and (3.2).

\vspace{0.5cm}
{\bf{Theorem 3.5.}}
\emph{For all odd integers $n\geq 1$ and real numbers $a\geq 1$, $x\in (0,\pi)$, we have
\begin{equation}
S^*_{n,a}(x)
\geq \sin(x).
\end{equation}
Equality holds if and only if $n=1$.}

\vspace{0.5cm}
\begin{proof}
Let $n\geq 1$ be odd and $a\geq 1$, $x\in (0,\pi)$.  Inequality (3.1) leads to
$$
2 S^*_{n,a}(x)
=S_{n,a}(x)+S_{n,a}(\pi-x)\geq \sin(x)+\sin(\pi-x)=2\sin(x).
$$
If equality holds in (3.7), then
$$
S_{n,a}(x)=\sin(x)
\quad\mbox{and}
\quad{S_{n,a}(\pi-x)=\sin(\pi-x).}
$$
From Theorem 3.1 we conclude that $n=1$.
\end{proof}

\vspace{0.5cm}
{\bf{Theorem 3.6.}}
\emph{For all even integers $n\geq 2$ and real numbers $a\geq 1$, $x\in (0,\pi)$, we have
\begin{equation}
S^*_{n,a}(x)
\geq 2 \sin(x).
\end{equation}
Equality holds if and only if $n=2$, $a=1$ or $n=4$, $a=1$, $x=\pi/2$.}

\vspace{0.5cm}
\begin{proof}
Let $n\geq 2$ be even and $a\geq 1$, $x\in (0,\pi)$. We apply (3.2) and obtain
$$
2 S^*_{n,a}(x)
=S_{n,a}(x)+S_{n,a}(\pi-x)
\geq 2\sin(x)(1+\cos(x))
+2\sin(\pi-x)(1+\cos(\pi-x))=4\sin(x).
$$
If $n=2$, $a=1$ or $n=4$, $a=1$, $x=\pi/2$, then equality is valid in (3.8). Conversely, if equality holds in (3.8), then
$$
S_{n,a}(x)=2\sin(x)(1+\cos(x))
\quad\mbox{and}
\quad{S_{n,a}(\pi-x)=2\sin(\pi-x)(1+\cos(\pi-x))
}.
$$
Applying Theorem 3.2 leads to $n=2$, $a=1$ or $n=4$, $a=1$, $x=\pi/2$.
\end{proof}

\vspace{0.5cm}
Now, we study the sine sum in two variables given in (1.4). The following two theorems offer improvements of inequality (1.2).

\vspace{0.5cm}
{\bf{Theorem 3.7.}}
\emph{For all odd integers $n\geq 1$ and real numbers $a\geq 1$, $x,y\in (0,\pi)$, we have
\begin{equation}
\Theta_{n,a}(x,y)\geq \sin(x)\sin(y).
\end{equation}
Equality holds if and only if $n=1$.}

\vspace{0.5cm}
\begin{proof}
Let $n\geq 1$ be  odd  and $a\geq 1$.
Since equality holds in (3.9) if $n=1$, we suppose that $n\geq 3$. Applying Theorem 3.1 gives for $x,y\in\mathbf{R}$ with $0<x-y<\pi$ and $0<x+y<\pi$:
$$
S_{n,a}(x-y)+
S_{n,a}(x+y)\geq \sin(x-y)+ \sin(x+y).
$$
This leads to
\begin{equation}
\sum_{j=1}^{n} {n+a-j\choose n-j} \sin(jx)\cos(jy)\geq \sin(x)\cos(y).
\end{equation}
Equality holds in (3.10) if and only if $n=3$, $a=1$ and  $x-y=2\pi/3$, $x+y=2\pi/3$,
that is,
$x=2\pi/3$, $y=0$. To obtain (3.9) we integrate both sides of (3.10) with respect to $y$. 

Let $x_0,y _0\in\mathbf{R}$ with
$0<y_0\leq x_0<\pi$. We consider two cases.

Case 1. $x_0\leq \pi/2$.\\
Let $0<y<y_0$. Then, $0<x_0-y<x_0+y<\pi$. It follows from (3.10) with $``>"$ instead of $``\geq"$ that
\begin{equation}
\int_0^{y_0}\sum_{j=1}^{n} {n+a-j\choose n-j} \sin(jx_0)\cos(jy)dy >\int_0^{y_0}
\sin(x_0)\cos(y)dy.
\end{equation}
This leads to (3.9) with $x=x_0$, $y=y_0$ and
 $``>"$ instead of $``\geq"$.
 
 Case 2. $\pi/2<x_0$.\\
 Case 2.1. $y_0\leq \pi-x_0$.\\
 Let $0<y<y_0$. Then, $0<x_0-y<x_0+y<\pi$, so that we obtain (3.11) and (3.9) 
 with $x=x_0$, $y=y_0$ and
 $``>"$ instead of $``\geq"$.\\
 Case 2.2. $\pi-x_0<y_0$.\\
 We set $x_1=\pi-x_0$ and $y_1=\pi-y_0$. Let $0<y<x_1$. Then, $0<y_1-y<y_1+y<\pi$. This leads to (3.11) with $x_1$ instead of $y_0$ and $y_1$ instead of $x_0$. Using
 $$
 \sin(jy_1)=(-1)^{j-1} \sin(j y_0)
 \quad\mbox{and}
 \quad{\sin(j x_1)=(-1)^{j-1} \sin(jx_0)}
 \quad{(j\in\mathbf{N})},
 $$
  we obtain (3.9) 
 with $x=x_0$, $y=y_0$ and
 $``>"$ instead of $``\geq"$.
 \end{proof}

\vspace{0.5cm}
{\bf{Theorem 3.8.}}
\emph{For all even integers $n\geq 2$ and real numbers $a\geq 1$, $x,y\in (0,\pi)$, we have
\begin{equation}
\Theta_{n,a}(x,y)\geq 2 \sin(x)\sin(y)(1+\cos(x) \cos(y)).
\end{equation}
Equality holds if and only if $n=2$, $a=1$.}

\vspace{0.5cm}
\begin{proof}
The proof is similar to that of Theorem 3.7. Therefore, we only offer a   proof sketch.
Since
$$
\Theta_{2,a}(x,y)=2\sin(x)\sin(y)
\Bigl(\frac{a+1}{2}+\cos(x)\cos(y)\Bigr),
$$
we conclude that if $n=2$, then (3.12) is valid with equality if and only if $a=1$.
Let $n\geq 4$. Using Theorem 3.2 we obtain for $x,y\in\mathbf{R}$ with
$0<x-y<\pi$ and $0<x+y<\pi$:
\begin{equation}
\sum_{j=1}^n {n+a-j\choose n-j}\sin(jx)\cos(jy)\geq 2\sin(x)\cos(y)+\sin(2x)\cos(2y).
\end{equation}
Equality holds in (3.13) if and only if $n=4$, $a=1$, $x=\pi/2$, $y=0$.

Let $x_0,y_0\in\mathbf{R}$ with $0<y_0\leq x_0<\pi$. We assume that $x_0\leq \pi/2$. If $0<y<y_0$, then $0<x_0-y<x_0-y<\pi$. Next, we integrate both sides of (3.13) (with $x=x_0$ and $``>"$ instead  of $``>"$) from $y=0$ to $y=y_0$. This gives
\begin{equation}
\Theta_{n,a}(x_0,y_0)>2\sin(x_0)\sin(y_0)(1+\cos(x_0) \cos(y_0)).
\end{equation}
By similar arguments we conclude that (3.14) also holds if $x_0>\pi/2$.
\end{proof}

\vspace{0.5cm}
 Let
$$
\Theta^*_{n,a}(x,y)=\sum_{j=1 \atop j \,  odd}^n
{n+a-j\choose  n-j}\frac{ \sin(jx)\sin(jy)}{j}.
$$
Applying Theorems 3.7 and 3.8
we obtain the following companions of (3.7) and (3.8).

\vspace{0.5cm}
{\bf{Theorem 3.9.}} \emph{For all real numbers $a\geq 1$ and $x,y\in (0,\pi)$ we have
\begin{equation}
\Theta^*_{n,a}(x,y)\geq \sin(x)\sin(y),
\quad\mbox{if}
\,\, {n} \,\, \mbox{is odd}
\end{equation}
and
\begin{equation}
\Theta^*_{n,a}(x,y)\geq 2 \sin(x)\sin(y),
\quad\mbox{if}
\,\, {n} \,\, \mbox{is even.}
\end{equation}
Equality holds in} (3.15) \emph{if and only if $n=1$ and in} (3.16) \emph{if and only if} $n=2$, $a=1$.

\vspace{0.5cm}
The proof of this theorem is quite similar to the proofs of Theorems 3.5 and 3.6, so that we omit the details.
We conclude this section with a generalization of inequality (1.2).

\vspace{0.5cm}
{\bf{Theorem 3.10.}}
\emph{Let $a$ be a  real number. The inequality
\begin{equation}
\Theta_{n,a}(x,y)>0
\end{equation}
holds for all integers $n\geq 1$ and real numbers $x,y\in (0,\pi)$ if and only if} $a\geq 1$.

\vspace{0.5cm}
\begin{proof}
Let $a\geq 1$. From Theorem 3.7
and 3.8 we conclude that $\Theta_{n,a}(x,y)$ is positive for all $n\geq 1$ and $x,y\in (0,\pi)$. Conversely, if (3.17) holds for all $n\geq 1$ and $x,y\in (0,\pi)$,
then we have
$$
\Theta_{2,a}(0,y)=0
\quad\mbox{and}
\quad{\frac{\partial}{\partial x} \Theta_{2,a}(x,y)\Big{|_{x=0}}=\sin(y)(1+a+2\cos(y))\geq 0}.
$$
It follows that
$$
1+a+2\cos(y)\geq 0.
$$
We let $y\rightarrow \pi$ and obtain $a-1\geq 0$.
\end{proof}

\vspace{0.5cm}
\section{Chebyshev polynomials}

The Chebyshev polynomials of the first and  second kind, $T_n(x)$ and $U_n(x)$ $(n=0,1,2,...)$, are polynomials in $x$ of degree $n$  defined by
\begin{equation}
T_n(x)=\cos(nt)
\quad\mbox{and}
\quad{
U_n(x)=\frac{\sin((n+1)t)}{\sin(t)},}
\quad{x=\cos(t),}
\quad{t\in[0,\pi].}
\end{equation}
They have remarkable applications in  numerical analysis, approximation theory and other branches of mathematics. The main properties of these functions can be found, for instance, in Mason and Handscomb \cite{MH}.

Using the notation (4.1) we obtain from (1.1) the inequality
\begin{equation}
\Lambda_{n,k}(x)=
\sum_{j=0}^n {n+k-j\choose n-j} U_j(x)>0
\quad{(0\leq n \in\mathbf{Z},
k\in \mathbf{N}; \, -1<x<1)}.
\end{equation}
The results presented in Section 3 lead to inequalities for Chebyshev polynomials of the
second kind. As examples we offer counterparts of Theorems 3.1 and 3.2 which provide refinements of (4.2).

\vspace{0.5cm}
{\bf{Theorem 4.1.}}
\emph{For all even integers $n\geq 0$ and real numbers $a\geq 1$, $x\in (-1,1)$, we have
$$
\Lambda_{n,a}(x)\geq 1.
$$
Equality holds if and only if $n=0$ or $n=2$, $a=1$, $x=-1/2$.}

\vspace{0.5cm}
{\bf{Theorem 4.2.}}
\emph{For all odd integers $n\geq 1$ and real numbers $a\geq 1$, $x\in (-1,1)$, we have
$$
\Lambda_{n,a}(x)\geq 2(1+x).
$$
Equality holds if and only if $n=1$, $a=1$ or $n=3$, $a=1$, $x=0$.}

\vspace{0.5cm}
The inequalities
\begin{equation}
0<\sum_{j=1}^n\frac{\sin(jx)}{j}<\pi-x
\quad{(n\in\mathbf{N}; 0<x<\pi)}
\end{equation}
are classical results in the theory of trigonometric polynomials. In 1910, Fej\'er conjectured the validity of the
left-hand side. The first proofs were published by Jackson \cite{J} in 1911 and  Gronwall \cite{G} in 1912. The right-hand side is due to Tur\'an \cite{T2}, who proved his inequality in 1938.
The following interesting companion of (4.3) was given
by Carslaw \cite{C} in 1917:
\begin{equation}
0<\sum_{j=0}^n \frac{\sin((2j+1)x)}{2j+1}\leq 1
\quad{(0\leq n\in\mathbf{Z}, 0<x<\pi)}.
\end{equation}
Equality holds if and only if $n=0$, $x=\pi/2$. See also \cite{AKO}.
We show that
  applications
of (4.3) and (4.4) lead to   integral inequalities for the Chebyshev polynomials.

 \vspace{0.5cm}
 \noindent
 {\bf{Theorem 4.3.}} \emph{For all integers $n\geq 1$ and real numbers $x\in (0,1)$, we have
 \begin{equation}
 \arccos(x)<\int_x^1 \frac{U_{2n}(t)}{\sqrt{1-t^2}}dt<\pi-\arccos(x).
 \end{equation}
 If $x\in (-1,0)$, then} (4.5) \emph{holds with $``>"$ instead of} $``<"$.
 
 \vspace{0.5cm}
 \begin{proof}
 Let $t\in (0,\pi)$ and
 $$
 F_n(t)=\sum_{j=1}^n \frac{\sin(2jt)}{j}.
 $$
 Then,
 $$
 F_n'(t)=2\sum_{j=1}^n \cos(2jt)
 =\frac{\sin((2n+1)t)}{\sin(t)}-1=U_{2n}(\cos(t))-1.
 $$
 This gives for $z\in (0,\pi)$:
 $$
 F_n(z)=\int_0^z F_n'(t) dt=\int_0^z U_{2n}(\cos(t))dt-z
 =\int_{\cos(z)}^1\frac{U_{2n}(t)}{\sqrt{1-t^2}}dt-z.
 $$
 If $x\in (0,1)$, then $\arccos(x)\in (0,\pi/2)$. Using 
 (4.3) leads to
 $$
 0<F_n(\arccos(x))=
 \int_{x}^1\frac{U_{2n}(t)}{\sqrt{1-t^2}}dt-\arccos(x)<\pi-2\arccos(x).
 $$
 This implies (4.5).
 
  We define for $x\in (-1,1)$:
 $$
 G_n(x)=\arccos(x)-
 \int_{x}^1\frac{U_{2n}(t)}{\sqrt{1-t^2}}dt
 \quad\mbox{and}
 \quad{w(x)=\arccos(x)-\frac{\pi}{2}.}
 $$
 Applying
 \begin{equation}
 w(x)+w(-x)=0
 \end{equation}
 and
 $$
 \int_{0}^1\frac{U_{2n}(t)}{\sqrt{1-t^2}}dt=\frac{\pi}{2}
 $$
 yields
 $$
 G_n(x)+G_n(-x)=\Bigl(2\int_0^1 -\int_x^1-\int_{-x}^1\Bigr)
 \frac{U_{2n}(t)}{\sqrt{1-t^2}}dt
 =\Bigl(\int_0^x-
 \int_{-x}^0\Bigr)
 \frac{U_{2n}(t)}{\sqrt{1-t^2}}dt.
$$
Since $U_{2n}$ is an even function, we obtain
$$
\int_{0}^x\frac{U_{2n}(t)}{\sqrt{1-t^2}}dt
=\int_{0}^x\frac{U_{2n}(-t)}{\sqrt{1-t^2}}dt
=\int_{-x}^0\frac{U_{2n}(t)}{\sqrt{1-t^2}}dt.
$$
Thus,
\begin{equation}
G_n(x)+G_n(-x)=0.
\end{equation}
Using (4.6) and (4.7) gives that the function
$$
H_n(x)=G_n(x)-2 w(x)
$$
satisfies
\begin{equation}
H_n(x)+H_n(-x)=0.
\end{equation}
From (4.5), (4.7) and (4.8) we conclude
 that for $x\in (-1,0)$ we have
$$
-G_n(x)=G_n(-x)< 0<H_n(-x)=-H_n(x).
$$
This leads to (4.5) with $``>"$ instead of $``<"$.
 \end{proof}

 \vspace{0.5cm}
{ \bf{Remark 4.4.}} The functions
$$
x\mapsto \int_x^1 \frac{U_n(t)}{\sqrt{1-t^2}}dt-\arccos(x)
\quad{(n=3,7,11)}
$$
 and
$${
x\mapsto \pi-\arccos(x)- \int_x^1 \frac{U_n(t)}{\sqrt{1-t^2}}dt}
\quad{(n=1,5,9)}
$$
 attain positive and negative values on $(0,1)$. This implies that Theorem 4.3 is in general not true for Chebyshev polynomials of odd degree.

\vspace{0.5cm}
We conclude this paper with a theorem which offers sharp upper and lower bounds for an integral involving the Chebyshev polynomials of the first and second kind.

\vspace{0.5cm}
{\bf{Theorem 4.5.}} \emph{For all integers $n\geq 0$ and real numbers $x\in (-1,1)$ we have
\begin{equation}
0<\int_x^1\frac{T_{n+1}(t) \, U_n(t)}{\sqrt{1-t^2}} dt\leq 1.
\end{equation}
Equality holds on the right-hand side if and only if} $n=0$, $x=0$.

\vspace{0.5cm}
\begin{proof}
Let $x\in (-1,1)$. Since $T_1(t)=t$ and $U_0(t)=1$, we obtain
$$
\int_x^1\frac{T_{1}(t) U_0(t)}{\sqrt{1-t^2}} dt
 =\sqrt{1-x^2}.
$$
This leads to (4.9) with $n=0$.

Let $n\geq 1$. We denote the sine sum in (4.4) by $I_n(x)$. Then,
 for $t\in (0,\pi)$,
 $$
 I_n'(t)=\sum_{j=0}^n\cos((2j+1)t)=\frac{\cos((n+1)t)\sin((n+1)t)}{\sin(t)}=T_{n+1}(\cos(t))\, U_n(\cos(t)).
 $$
 Hence, we get for $z\in (0,\pi)$:
 $$
 I_n(z)=\int_0^z I_n'(t)dt=\int_{\cos(z)}^1 \frac{T_{n+1}(t)\, U_n(t)}{\sqrt{1-t^2}} dt,
 $$
 so that (4.4) (with ``$<$" instead of ``$\leq $") gives
 $$
 0<I_n(\arccos(x))=\int_x^1 \frac{T_{n+1}(t)\,U_n(t)}{\sqrt{1-t^2}} dt <1.
 $$
 This completes the proof.
\end{proof}

\vspace{0.5cm}
{\bf{Acknowledgement.}} We thank the referee for the careful reading of the paper and for  helpful comments.

\end{document}